\documentclass[12pt]{amsart}
\usepackage{amsmath}
\usepackage{enumerate}
\usepackage[pdftex]{graphicx}
\usepackage{cite}
\usepackage{hyperref}
\theoremstyle{remark}

\newcommand{\rref}{\ref}                 
\newtheorem{Theorem}{Theorem}[section]
\newtheorem{Corollary}[Theorem]{Corollary}

\newtheorem{Conjecture}[Theorem]{Conjecture}

\def\eqref#1{(\rref{#1})}                

\long\def\commentout#1{}

\begin{document}


\title{Single-colored ADO-3 invariant is a specialization of Links-Gould invariant for closures of 5-braids}


\author{Nurdin Takenov}
\address{Department of Mathematics, 303 Lockett Hall, Louisiana State University, Baton Rouge, LA 70803}
\email{tnurdi1@lsu.edu}
\urladdr{www.math.lsu.edu/$\sim$takenov} 

\begin{abstract}
We prove that for knots and links that are closures of 5-braids, single-colored ADO-invariant  of third order coincides with Links-Gould polynomial with specific choice of variables. We also conjecture that this is true for any knots and links and provide some arguments for that.
\end{abstract}

\maketitle

\tableofcontents

\section{Introduction}

In this article we study the relation between two knot invariants, namely ADO polynomial(also known as colored Alexander polynomial) and Links-Gould polynomial. ADO polynomials $ADO_N(L;e^{\lambda_1}, \dots)$ of a link $L$ is a family of quantum invariants, indexed by positive integers $N=2,3,\dots$\cite{akutsu92}, \cite{murakami08}. In ADO polynomials each connected components is assigned separate color, however we will consider the case when all components colored by the same color $\lambda$. In this article, as we explain later, we will use variable $t=q^{\lambda}$. We will still call this single-colored version $ADO_N(L;t)$, since in the most important case of knots there is no difference. 

Links-Gould polynomial $LG(L; t_0, t_1)$ is another quantum invariant \cite{links92}, \cite{dewit99},\cite{ishii04}\footnote{We will use Ishii's notation.}. The main result of this paper is a relation between ADO polynomial for $N=3$ and Links-Gould polynomial:

\begin{Theorem}\label{maintheorem} If link $L$ is a closure of a 5-braid, then 
\begin{align}
ADO_3(L;t) = LG(L; t^2,\omega^2 t^{-2} ),
\end{align}

where 
\begin{align*}
\omega = e^{\pi i/3}=\frac{1+i\sqrt{3}}{2}.
\end{align*}
\end{Theorem}

As a corollary we get the following property about values of ADO polynomial at $t=1$:

\begin{Corollary}\label{maincorollary}
If link $L$ is a closure of a 5-braid, then
\begin{align}
ADO_3(L;1) & = ADO_3(L;\omega) =\begin{cases}1, \text{ if }L \text{ is a knot,}\\
0, \text{ otherwise}.
\end{cases}
\end{align}
\end{Corollary}

We conjecture that these statements hold in general, for any links:

\begin{Conjecture}\label{mainconjecture} For any link $L$ we have
\begin{align}
ADO_3(L;t) & = LG(L; t^2,\omega^2 t^{-2} ),\\
ADO_3(L;1) & = ADO_3(L;\omega)  =\begin{cases}1, \text{ if }L \text{ is a knot,}\\
0, \text{ otherwise}.
\end{cases}
\end{align}
\end{Conjecture}

Similar relations between Links-Gould polynomial and regular Alexander polynomial was proved before\cite{kohli17}, so it's natural to assume that there are more identities to be proven.

\section{Knot invariants} Both ADO polynomial and Links-Gould polynomials are operator invariants(general overview can be seen in \cite{ohtsuki02}, Chapters 3 and 4) of oriented links. To construct invariants we need to specify vector space $V$ and then define two maps, one corresponding to crossing move, given by $R$-matrix, $R:V\otimes V\to V\otimes V$ and map $h: V\to V$, which defines cup and cap operators. 

We represent link $L$ as a $(1,1)$-tangle(Usually in simiilar constructions links represented as a $(0,0)$-tangle, but we avoiding this since for $(0,0)$-tangles our invariants will nullify). The easiest way to do is to consider knot as closure of a braid $\beta$ and then cut the rightmost strand.
\begin{figure}[h]
\includegraphics[scale=0.5]{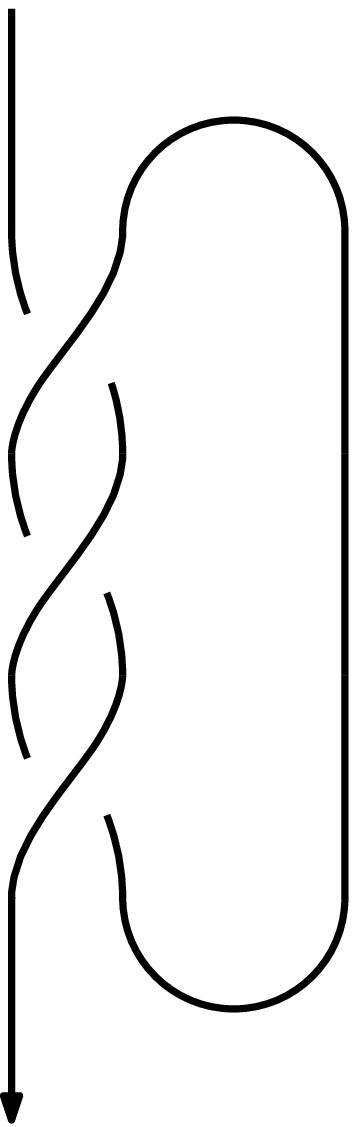}
\caption{Trefoil represented as $(1,1)$ tangle.}
\end{figure}

Then we slice the diagram and consider the operator invariant resulting from applying crossing moves, cup and cap operators.

\begin{figure}[h]
\includegraphics[scale=1.0]{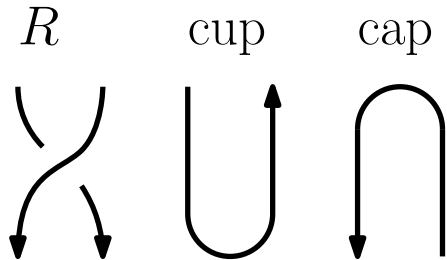}
\caption{$R$, cup and cap operators.}
\end{figure}
We will be reading sliced diagrams from bottom to top, therefore cup operator is a map from ground field to $V\otimes V^{*}$. In the end we will get a linear operator $O_{\beta}:V\to V$, which should be scalar multiplication. The value of that scalar is a link invariant, given proper choice of $R$ and $h$. 

Another way to express it would be to note that, given proper map $R$, we can define representation $\phi^n_R$ of braid group $B_n$ by defining 
\begin{align*}
\phi^n_R(s_i) = \text{id}_V^{\otimes i-1} \otimes R \otimes \text{id}_V^{\otimes n-i-1},
\end{align*}
where $s_i$ - Artin generators of braid group, $i=1,\dots, n-1$.
 $O_{\beta}$ is then a partial trace: 
 \begin{align*}
 O_{\beta} = \text{trace}_{2,\dots, n}\left(\left(\text{id}_V\otimes h^{\otimes(n-1)} \right)\circ \phi^n_R(\beta)\right): V\to V
 \end{align*}
 
 Then there exist scalar $c_{\beta}$, such that $O_{\beta} = c\text{ id}_V$. Given the proper choice of $R$ and $h$, the value $c_{\beta}$ will not change under Markov moves, so then $c_{\beta}$ would be knot invariant.

We will provide brief description of these invariants.

\subsection{ADO polynomial}

ADO polynomial was introduced in 1992 by Akutsu, Deguchi and Ohtsuki\cite{akutsu92}. We will however use the description given by Murakami in 2008 \cite{murakami08} with some modifications and restricting it to the case when all components colored byt the same color $\lambda$. Modifications required in order to make the ADO polynomial invariant of framing, also we use different normalization. Moreover, we will focus our attention on case $N=3$, the corresponding invariant we are going to refer as ADO-3 polynomial. 

First, we will construction as is given by Murakami. For a given integer $N$, such that $N\ge 2$, we fix $q = \exp(\pi i/N)$ to be a primitive $2N$-th root of unity and also introduce variable $\lambda$, color \footnote{In case of a link, we assign the same color $\lambda$ to all components.}. Then we consider an $N$-dimensional vector space $V_{\lambda}$ with a basis $\{v_0, v_1,\dots, v_{N-1}\}$. Action of $R$-matrix for ADO polynomial is given by:
\begin{align}
& R(v_i\otimes v_j)  = \\
& \sum_n q^{\frac{1}{2}(\lambda -2i-2n)(\lambda -2j+2n)+n(n-1)/2} \frac{\{i+n;n\}\{\lambda -j+n; n\}}{\{n;n\}} v_{j-n}\otimes v_{i+n} \nonumber
\end{align}
where
\begin{align*}
\{a\} & = q^a-q^{-a}, \\
\{x;n\} & =\begin{cases}
\prod^{n-1}_{i=0} \{x-i\}, n>0\\
1, n=0.
\end{cases}
\end{align*}

Cup and cap operators are given as:
\begin{align*}
\text{cup}(1) & = \sum_{i=0}^{N-1} q^{(N-1)\lambda+2i} v_i\otimes v_i^{*}\\
\text{cap}(v_i\otimes v_j^{*}) & = \delta_{ij}
\end{align*}

Then if link $L$ is a closure of a $(1,1)$ tangle $T$, and $O_T^{N}(\lambda)$ is the corresponding operator, ADO-invariant is defined by 
\begin{align*}
ADO_N(L;\lambda) \text{Id}_V= \frac{q^{-\frac{\lambda(\lambda+2-2N)}{2}f}}{\{\lambda+N; N-1\}} O_T^{N}(\lambda),
\end{align*}
where $f$ is framing, number of positive crossings minus number of negative crossings(Theorem 5 and Proposition 6 from \cite{murakami08}). 

We will make the following modifications to this construction: first, we remove the normalization factor $\{\lambda+N; N-1\}^{-1}$, since we want our invariant to be $1$ in case of unknot. Second, to remove framing correction factor, we will multiply $R$ matrix by the factor $q^{-\frac{\lambda(\lambda+2-2N)}{2}}$, this way framing correction will happen automatically. These modifications are similar to the ones made in \cite{ito16}, but not identical.

The modified version looks like this: action of $R$-matrix for ADO polynomial is given by:
\begin{align}
& R(v_i\otimes v_j)  = \\
& q^{(N-1)\lambda - (i+j)\lambda}\sum_n q^{2(i+n)(j-n)+n(n-1)/2} \frac{\{i+n;n\}\{\lambda -j+n; n\}}{\{n;n\}} v_{j-n}\otimes v_{i+n} \nonumber.
\end{align}

Cup and cap operators are given as:
\begin{align*}
\text{cup}(1) & = \sum_{i=0}^{N-1} q^{(N-1)\lambda+2i} v_i\otimes v_i^{*}\\
\text{cap}(v_i\otimes v_j^{*}) & = \delta_{ij}
\end{align*}

Then if link $L$ is a closure of a $(1,1)$ tangle $T$, and $O_T^{N}(\lambda)$ is the corresponding operator, ADO-invariant is defined by 
\begin{align*}
ADO_N(L;\lambda) \text{Id}_V=  O_T^{N}(\lambda).
\end{align*}

This is the version that we will be using from now on. $ADO_N(L;\lambda)$ is Laurent polynomial for $q$ and $q^{\lambda}$, so we will work instead with the variable $t = q^{\lambda}$. We can note that $ADO_N(L;\lambda)$ takes values in $\mathbb{Z}[q, t^{\pm 1}]$.

In particular, when $N=3$, $R$-matrix for ADO-polynomial is 
\begin{align*}
R = \left(
\begin{array}{ccccccccc}
 t^2 & 0 & 0 & 0 & 0 & 0 & 0 & 0 & 0 \\
 0 & t^2-1 & 0 & t & 0 & 0 & 0 & 0 & 0 \\
 0 & 0 & \frac{\left(t^2-1\right) \left(2 t^2-i \sqrt{3}+1\right)}{2 t^2} & 0 & t-\frac{1}{t} & 0 & 1 & 0 & 0 \\
 0 & t & 0 & 0 & 0 & 0 & 0 & 0 & 0 \\
 0 & 0 & \frac{1}{t}+\frac{1}{2} \left(1+i \sqrt{3}\right) t & 0 & \frac{-1+i\sqrt{3}}{2} & 0 & 0 & 0 & 0 \\
 0 & 0 & 0 & 0 & 0 & -\frac{2 t^2-i \sqrt{3}+1}{2 t^2} & 0 & -\frac{1+i \sqrt{3}}{2 t} & 0 \\
 0 & 0 & 1 & 0 & 0 & 0 & 0 & 0 & 0 \\
 0 & 0 & 0 & 0 & 0 & -\frac{1+i \sqrt{3}}{2 t} & 0 & 0 & 0 \\
 0 & 0 & 0 & 0 & 0 & 0 & 0 & 0 & \frac{-1+i\sqrt{3}}{2 t^2} \\
\end{array}
\right)
\end{align*}

Here $q = \omega = \exp(\pi i/3)$, $t = q^{\lambda}$, the basis is ordered as 
\begin{align*}
v_0\otimes v_0,v_0\otimes v_1, v_0\otimes v_2, v_1\otimes v_0,v_1\otimes v_1, v_1\otimes v_2, v_2\otimes v_0,v_2\otimes v_1, v_2\otimes v_2. 
\end{align*}

\subsection{Links-Gould polynomial}

First Links-Gould invariant $LG(L; t_0, t_1)$ was introduced in 1992 by Links and Gould\cite{links92} and more thoroughly studied by De Wit, Kauffman and Links in 1999\cite{dewit99}. It is the first in a series of Links-Gould invariants $LG^{n,m}$ and denoted as $LG^{2,1}$ to distinguish it from others. For example $LG^{1,1}$ is Alexander polynomial\cite{viro07}. We will however only work with the first Links-Gould invariant, so we will refer to it just as $LG$. We can also note that Links-Gould polynomial takes values in $\mathbb{Z}[t_0^{\pm 1}, t_1^{\pm 1}]$ (Theorem 1, \cite{ishii06}).

In case of a Lanks-Gould polynomial, the vector space $W$ is 4-dimensional, with basis $\{e_0, e_1, e_2, e_3\}$. The map $h$ is given by matrix 
\begin{align*}
\left(
\begin{array}{cccc}
 \frac{1}{t_0} & 0 & 0 & 0 \\
 0 & -t_1 & 0 & 0 \\
 0 & 0 & -\frac{1}{t_0} & 0 \\
 0 & 0 & 0 & t_1 \\
\end{array}
\right),
\end{align*}

To write the $R$-matrix for Links-Gould polynomial we denote $Y = \sqrt{\left(t_0-1\right) \left(1-t_1\right)}$. Then $R$-matrix for Links-Gould polynomial is 

\begin{align*}
\setlength{\arraycolsep}{.5\arraycolsep}
\begin{pmatrix}
 t_0 & 0 & 0 & 0 & 0 & 0 & 0 & 0 & 0 & 0 & 0 & 0 & 0 & 0 & 0 & 0 \\
 0 & 0 & 0 & 0 & \sqrt{t_0} & 0 & 0 & 0 & 0 & 0 & 0 & 0 & 0 & 0 & 0 & 0 \\
 0 & 0 & 0 & 0 & 0 & 0 & 0 & 0 & \sqrt{t_0} & 0 & 0 & 0 & 0 & 0 & 0 & 0 \\
 0 & 0 & 0 & 0 & 0 & 0 & 0 & 0 & 0 & 0 & 0 & 0 & 1 & 0 & 0 & 0 \\
 0 & \sqrt{t_0} & 0 & 0 & t_0-1 & 0 & 0 & 0 & 0 & 0 & 0 & 0 & 0 & 0 & 0 & 0 \\
 0 & 0 & 0 & 0 & 0 & -1 & 0 & 0 & 0 & 0 & 0 & 0 & 0 & 0 & 0 & 0 \\
 0 & 0 & 0 & 0 & 0 & 0 & t_0 t_1-1 & 0 & 0 & -\sqrt{t_0 t_1} & 0 & 0 & -\sqrt{t_0 t_1} Y & 0 & 0 & 0 \\
 0 & 0 & 0 & 0 & 0 & 0 & 0 & 0 & 0 & 0 & 0 & 0 & 0 & \sqrt{t_1} & 0 & 0 \\
 0 & 0 & \sqrt{t_0} & 0 & 0 & 0 & 0 & 0 & t_0-1 & 0 & 0 & 0 & 0 & 0 & 0 & 0 \\
 0 & 0 & 0 & 0 & 0 & 0 & -\sqrt{t_0 t_1} & 0 & 0 & 0 & 0 & 0 & Y & 0 & 0 & 0 \\
 0 & 0 & 0 & 0 & 0 & 0 & 0 & 0 & 0 & 0 & -1 & 0 & 0 & 0 & 0 & 0 \\
 0 & 0 & 0 & 0 & 0 & 0 & 0 & 0 & 0 & 0 & 0 & 0 & 0 & 0 & \sqrt{t_1} & 0 \\
 0 & 0 & 0 & 1 & 0 & 0 & -\sqrt{t_0 t_1} Y & 0 & 0 & Y & 0 & 0 & Y^2 & 0 & 0 & 0 \\
 0 & 0 & 0 & 0 & 0 & 0 & 0 & \sqrt{t_1} & 0 & 0 & 0 & 0 & 0 & t_1-1 & 0 & 0 \\
 0 & 0 & 0 & 0 & 0 & 0 & 0 & 0 & 0 & 0 & 0 & \sqrt{t_1} & 0 & 0 & t_1-1 & 0 \\
 0 & 0 & 0 & 0 & 0 & 0 & 0 & 0 & 0 & 0 & 0 & 0 & 0 & 0 & 0 & t_1 \\
\end{pmatrix}
\end{align*}

Then the Links-Gould polynomial is calculated according to general description given above for operator invariants.

\section{Common properties of ADO and Links-Gould invariants}
\label{common}

In this section we will list common properties of ADO and Links-Gould invariants, some of which we will use to prove the main result. 
\subsection{Cubic skein relation}

By direct computation, it can be checked that $R_{ADO-3}$, $R$-matrix for ADO polynomial when $N=3$ satisfied the cubic relation 
\begin{align*}
R^3 = -\omega ^2\text{Id}_9+R^2 \left(\frac{\omega ^2}{t^2}-1+t^2\right)+R \left(\frac{\omega ^2}{t^2}-\omega ^2+t^2\right).
\end{align*}

In fact this relation follows from a minimal polynomial for this particular matrix. This relation for $R$-matrices gives us cubic skein relation:
\begin{align}
\label{adoskein}
ADO_3(L_3;t) & = -\omega^2 ADO_3(L_0;t)+\left(\frac{\omega^2}{t^2}-\omega^2+t^2\right)ADO_3(L_1;t)+\\
\nonumber & + \left(\frac{\omega^2}{t^2}-1+ t^2\right)ADO_3(L_2;t)
\end{align}

where $L_k$ are links that look identically except for some ball, inside which two strands make $k$ twists. Pictorially it looks like this: 
\begin{figure}[h]
\includegraphics[scale=0.8]{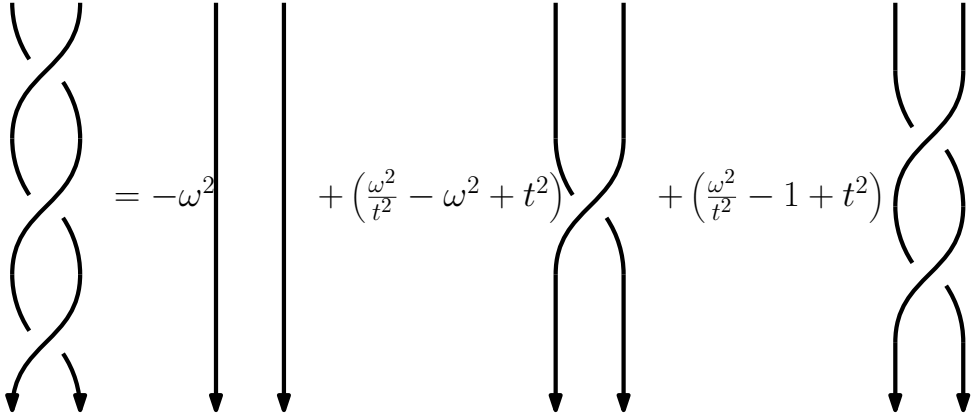}
\caption{Cubic skein relation for ADO-3 invariant.}
\end{figure}

Links-Gould polynomial famously satisfies skein relation
\begin{align}
\label{lgskein}
LG(L_3) & + (1-t_0-t_1) LG(L_2) + \\
\nonumber & + (t_0 t_1-t_0-t_1) LG(L_1) + t_0 t_1 LG(L_0)=0
\end{align}

This relation is given in \cite{dewit99}, p.170, but we use different set of variables, the one is used in \cite{ishii2004algebraic}. The crucial fact is that skein relations \ref{adoskein} and \ref{lgskein} become the same if we specify variables in Links-Gould variable to $t_0=t^2$, $t_1 = \omega^2 t^{-2}$. Therefore ADO polynomial $ADO_3(L;t)$ and specialized Links-Gould polynomial $LG(L; t^2,\omega^2 t^{-2} )$ satisfy the same cubic skein relation, the one shown on the picture above.

\subsection{Three strand relation}

This property is not used in the proof, but it may be useful to proving the main conjecture. Ishii in \cite{ishii04}, Proposition 1, described several three-strand relations for Links-Gould polynomial. Let's consider one of them, namely relation (3.4) in Proposition 1. Ishii introduced new operators $Q_i$ given by 
\begin{align*}
Q_0 & = \frac{t_0}{t_0-t_1}R +\frac{t_0(1-t_1)}{t_0-t_1} \text{id}_{V\otimes V} - \frac{t_0 t_1}{t_0-t_1} R^{-1},\\
Q_1 & = \frac{t_1}{t_1-t_0}R +\frac{t_1(1-t_0)}{t_1-t_0} \text{id}_{V\otimes V} - \frac{t_0 t_1}{t_1-t_0} R^{-1}.
\end{align*} 

Then the following relation holds:
\begin{align*}
(Q_0\otimes \text{id}_V)(\text{id}_V\otimes Q_1)(Q_1\otimes \text{id}_V) = \frac{t_1(t_0-1)}{t_0 \left(1-t_1\right)}(Q_0\otimes \text{id}_V)(\text{id}_V\otimes Q_0)(Q_1\otimes \text{id}_V)
\end{align*}

By direct calculation, it can be shown that if we specialize this relation, namely put $t_0=t^2$, $t_1 = \omega^2 t^{-2}$, the same relation holds for ADO polynomial for $N=3$. The same is true for other relations given by Ishii in Proposition 1. These calculations are done on computer. Details are given in the Appendix. 

\subsection{Invariants vanish on split links}

Another common property of a ADO and Links-Gould invariants is that they both vanish on split links. This follows from the fact that operator invariant for $(0,0)$ tangle is zero for both ADO and Links-Gould construction(the result for vanishing invariant of $(0,0)$-tangle can be found in \cite{akutsu92}, Proposition 4.3 and \cite{dewit00}, section 3.5). The exact argument for LG invariant can be found in \cite{dewit00}, section 3.5, but the same argumentation works for ADO-invariant.

\section{Proof of the main theorem}

We are considering the knots that are closures of 5-braids. Since both ADO and specialized Links-Gould polynomials satisfy the same cubic skein relation, we only need to consider the quotient of group of 5-braids by the cubic relation. This quotient is finite\cite{marin12}, so we can check the values of ADO and LG invariant on all of them. This part is done by computer calculations. 

\subsection{Quotient of braid group algebra $RB_5$ by cubic relation} 
\label{sn}This section closely follows Marin's paper \cite{marin12} and Marin and Wagner's paper \cite{marin13}. Let us consider braid group algebra $RB_n$, where $R = \mathbb{Z}[\omega^2, \omega^2/ t^2+ t^2] = \mathbb{Z}[\omega, \omega^2/ t^2+ t^2]$. Then we can consider $A_n$, quotient of $RB_n$ by cubic relation  
\begin{align}
\label{cubic}
s_i^3 - \left(\frac{\omega^2}{t^2}-1+ t^2\right)s_i^2 - \left(\frac{\omega^2}{t^2}-\omega^2+t^2\right)s_i + \omega^2 = 0,
\end{align}

where $s_i$ are Artin generators of a braid group. We want to use here results from \cite{marin12}, namely theorems 1.2 and 4.1. The problem is that in our case, the coefficients $a,b,c$(in Marin's notation) of a cubic relation are not independent, for example there are relations $a-b+c+1=0$, $c^3=1$. Still, most of the arguments can be adapted to our case, the only major difference would be that $A_n$, instead of a being free $R$-module, will be quotient of a free $R$-module generated by elements described in theorems 1.2 and 4.1 of \cite{marin12}. If we define this generating set as $S_n$, we have:
\begin{align*}
S_2 & = \{1, s_1, s_1^{-1}\},\\
S_3 & = S_2 \sqcup S_2 s_2^{\pm 1}S_2 \sqcup S_2 s_2^{-1}s_1s_2^{-1}
\end{align*}

Using proposition 4.8 of \cite{marin12}, we also get the description of $S_4$:
\begin{align*}
S_4 & = U S_3, \text{ where }\\
U & = \{1,s_3^{-1} s_2 s_1^{-1} s_2 s_3^{-1}, s_3 s_2^{-1} s_1 s_2^{-1} s_3, s_3^{\pm 1}, s_3^{\pm 1} s_2^{\pm 1},s_3^{\pm 1} s_2^{\pm 1} s_1^{\pm 1},\\
& s_3^{\pm 1} s_2^{-1} s_1 s_2^{-1},s_3 s_2^{-1} s_3, s_3 s_2^{- 1} s_3 s_1^{\pm 1}, s_3 s_2^{- 1} s_3 s_1 s_2^{- 1}s_1, s_3 s_2^{- 1} s_3 s_1^{\pm 1} s_2^{\pm 1}\}
\end{align*}

\subsection{Extending knot invariants to $RB_n$ and $A_n$}

Any knot invariant can be defined as a function on $B_n$ by assigning to an element of braid group the value of a knot invariant of a braid closure. Then, it can be linearly extended to $RB_n$. Since both $ADO_3(L;t)$ and $LG(L; t^2, w^2 t^{-2})$ invariants satisfy cubic relation \ref{cubic}, it means that they well-defined on $A_n$. Abusing the notation, we will use the same notation for extended functions on $A_n$. We want to show that these link invariants coincide on $B_5$. In light of the previous, equivalently we can try to prove that these link invariants coincide on $A_5$. But, since $A_5$ is finitely generated module over $R$, we only need to show that $ADO_3(L;t)$ and $LG(L; t^2, w^2 t^{-2})$ invariants coincide on $S_5$, finite generating set for $A_5$. We will do things slightly differently: we will prove the statement for $S_4$(and consequently for $A_4$ and $B_4$) and then we will use this to prove the result for $A_5$.

\subsection{Checking the statement of main theorem for $A_5$}
\label{a5}

Now we need to check the statement of main theorem for $A_5$. We will try to simplify calculations, for that we need first to take a look at $S_4$. Using the description of $S_4$ given above, we can generate all elements of $S_4$ and check the statement of conjecture for them. The corresponding calculations are performed on a computer, the details are given in an Appendix.

A description of $A_5$ is given in a theorem 6.21 from \cite{marin12}, we have:
\begin{align*}
A_5 & = A_4 + A_4 s_4 A_4 + A_4 s^{-1}_4 A_4 + A_4 s_4 s^{-1}_3 s_4 A_4 + A_4 s^{-1}_4 s_3 s^{-1}_2 s_3 s^{-1}_4 A_4 +\\
& + A_4 s_4 s^{-1}_3 s_2 s^{-1}_3 s_4 A_4 + A_4 s^{-1}_4 w^+ s^{-1}_4 A_4 + A_4 s_4 w^- s_4 A_4 + A_4s^{-1}_4 w^- s^{-1}_4 A_4 +\\
& + A_4 s_4 w^+ s_4 A_4 + s_4 w^- s_4 w^- s_4 A_4 + s_4 w^+ s^{-1}_4 w^+ s_4 A_4 + s^{-1}_4 w^- s_4 w^- s^{-1}_4 A_4,
\end{align*}

where $w^+ = s_3 s^{-1}_2 s_1 s^{-1}_2 s_3$, $w^- = s^{-1}_3 s_2 s^{-1}_1 s_2 s^{-1}_3$. If can prove the main conjecture for separate subsets(like $A_4$, $A_4 s_4 A_4$, etc) then this proves the main conjecture for the whole of $A_5$. Let's consider subsets separately.

For $A_4\subset A_5$: we consider here the inclusion of $A_4$ into $A_5$, which basically means adding one additional strand, unconnected to others. Then, each element of $A_4$ is a linear combination of $5$-braids with one strand not interacting with others. If we take closure of such braid, we would get a split link, for which both link invariants are zero. Therefore for all elements of $A_4\subset A_5$ our invariants coincide.

For elements of the type $A_4 s_4^{\pm 1} A_4$ we can use second Markov move to reduce them to elements of $A_4$.

\begin{figure}[h]
\includegraphics[scale=1.0]{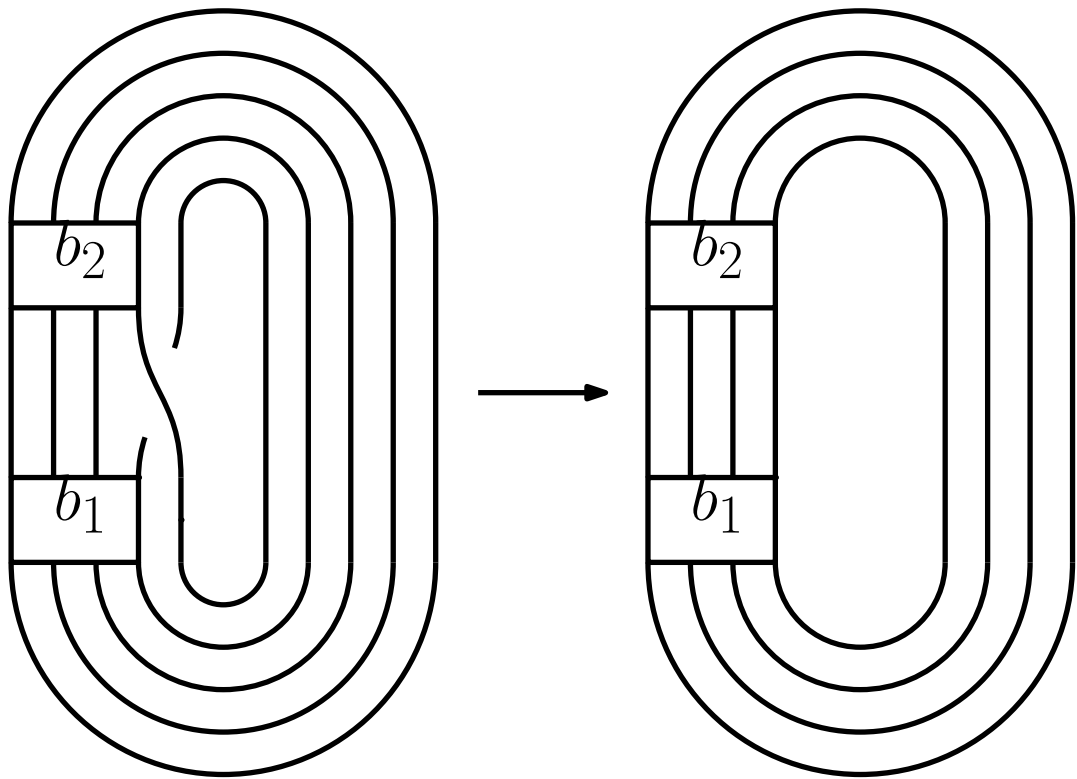}
\caption{Second Markov move.}
\end{figure}

Namely, every element of the type $A_4 s_4^{\pm 1} A_4$ can be written as a linear combination of elements of the form $b_1 s_4^{\pm 1} b_2$, where $b_1$ and $b_2$ are some 4-braids. When taking closure, using the second Markov move we can remove the middle term $s_4^{\pm 1}$, and we are left with a closure of a 4-braid $b_1 b_2$. But we already proved the statement of the theorem for closures of 4-braids, so it proves the statement for elements of the type $A_4 s_4^{\pm 1} A_4$.

Now let us consider the elements of the type $A_4 w A_4$, where $w$ is a 5-braid. Using the first Markov move we can show that they can be reduced to elements of the type $A_4 w$. Namely, every element of the type $A_4 w A_4$ can be written as linear combination of the elemenst of the form $b_1 w b_2$, where $b_1$, $b_2$ - some 4-braids. We need to show that ADO-3 and specialized Links-Gould polynomial coincide for the closure of a braid $b_1 w b_2$. Applying the first Markov move, we note that the closure of a braid $b_1 w b_2$ is the same as a closure of a braid $b_2 b_1 w$. But $b_2 b_1 w \in A_4 w$. So it is enough to show that ADO-3 and specialized Links-Gould polynomials coincide on $A_4 w$.

\begin{figure}[h]
\includegraphics[scale=1.0]{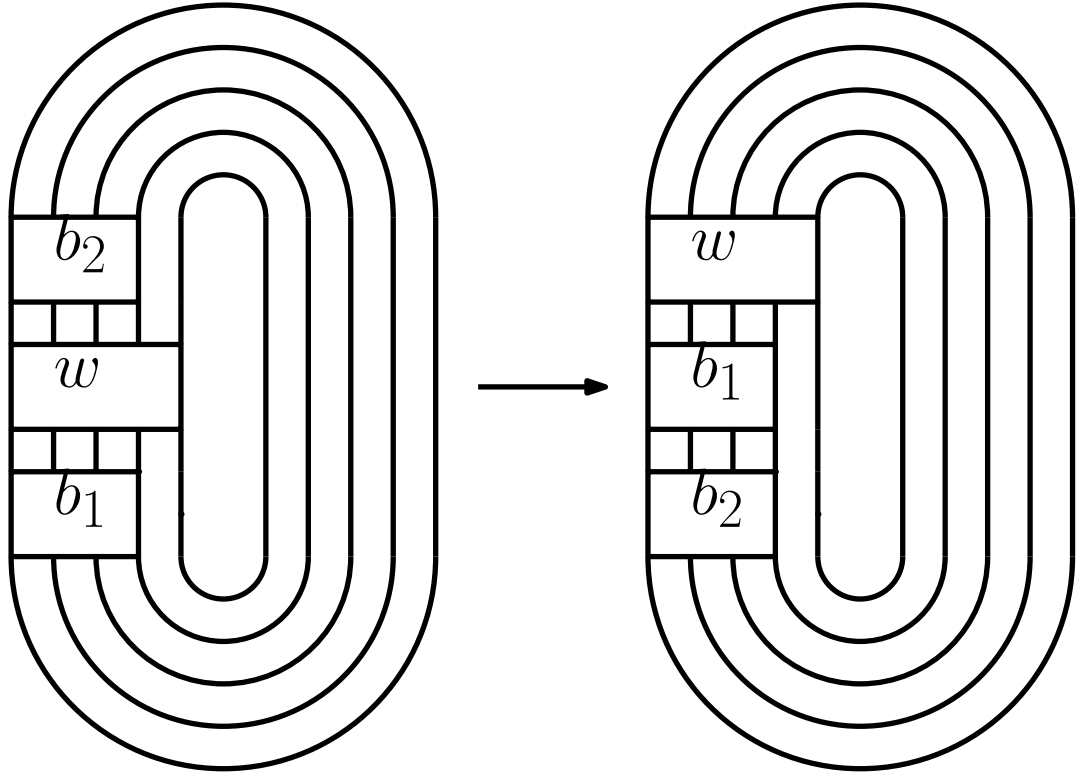}
\caption{First Markov move.}
\end{figure}

Since $A_4 w$ as an $R$-module is generated by $S_4 w$, it is enough to check the ADO-3 and specialized Links-Gould polynomials for for the set $S_4 w$. Given the description of $A_5$ given above, it means that in order to prove the statement of the theorem for $A_5$ it is enough to prove it for the sets:
\begin{align*}
& S_4 s_4 s^{-1}_3 s_4, S_4 s^{-1}_4 s_3 s^{-1}_2 s_3 s^{-1}_4, S_4 s_4 s^{-1}_3 s_2 s^{-1}_3 s_4, S_4 s^{-1}_4 w^+ s^{-1}_4, S_4 s_4 w^- s_4, \\
& S_4 s^{-1}_4 w^- s^{-1}_4 , S_4 s_4 w^+ s_4 , s_4 w^- s_4 w^- s_4 S_4, s_4 w^+ s^{-1}_4 w^+ s_4 S_4 ,  s^{-1}_4 w^- s_4 w^- s^{-1}_4 S_4,
\end{align*}
where $w^+ = s_3 s^{-1}_2 s_1 s^{-1}_2 s_3$, $w^- = s^{-1}_3 s_2 s^{-1}_1 s_2 s^{-1}_3$.

The calculations are done on a computer. The details are given in the Appendix.

\section{Proof of the corollary}

The corollary is a direct consequence of properties of a Links-Gould polynomial. According to Theorem 7 from \cite{ishii06}, which states that for any link $L$:
\begin{align*}
LG(L; t_0,1) = LG(L; 1,t_1) = \begin{cases}1, \text{ if }L \text{ is a knot,}\\
0, \text{ otherwise}.
\end{cases}
\end{align*}

Consequently we have:
\begin{align*}
ADO_3\left(L;1\right) & = LG\left(L; 1, \omega^2 \right) = \begin{cases}1, \text{ if }L \text{ is a knot,}\\
0, \text{ otherwise},
\end{cases}\\
ADO_3\left(L;\omega\right) & = LG\left(L; \omega^2, 1\right) = \begin{cases}1, \text{ if }L \text{ is a knot,}\\
0, \text{ otherwise}.
\end{cases}\\
\end{align*}

\section{Possible ways to prove conjecture}

The methods described here can be extended to closures of $5$-braids, however extension to six or more strands is difficult, since $A_n$ doesn't have finite basis for $n>5$. Another way to prove conjecture in full extent would be to properly study the $R$-matrices for Links-Gould and ADO polynomials, similar to \cite{kohli17}. 

One of the most promising ways is to use the main result of \cite{marin13}. In this paper Marin and Wagner showed that Links-Gould polynomial is completely defined by 3 skein-like relations, $r_1$, $r_2$ and $r_3$. The first relation, $r_1$ is a cubic relation \ref{lgskein}. As we have seen $ADO_3$ polynomial satisfies a specialization of this cubic relation, \ref{adoskein}.

The second relation, $r_2$ is one of the relations (3.4)-(3.6) from \cite{ishii04}. As we mentioned in a subsection about three strand relations, by direct computation in Mathematica we can show that $ADO_3$ satisfied the same relation with the same specialization. 

If we could show that $ADO_3$ satisfies the specialized relation $r_3$ this would prove the conjecture in full extent. Unfortunately, the explicit form of relation $r_3$ is still not known.

Another fact that implies validity of the main conjecture is the fact that ADO invariant is almost-symmetric:
\begin{align*}
ADO_3(L;t) & = ADO_3\left(L;\omega t^{-1}\right),\\
\end{align*}
It was proved recently by Martel and Willetts \cite{martel21} in more general form, but it also can be shown to be consequence of a main conjecture and a properties of a Links-Gould invariant:
\begin{align*}
ADO_3\left(L;\omega t^{-1}\right) & =  LG\left(L; \left(\omega t^{-1}\right)^2,\omega^2 \left(\omega t^{-1}\right)^{-2} \right) =\\
& = LG\left(L; \omega^2 t^{-2}, t^2 \right) = LG\left(L; t^2, \omega^2 t^{-2}\right) =\\
& = ADO_3\left(L;t\right)
\end{align*}

\bibliography{ado3}{}
\bibliographystyle{alpha}

\appendix
\section{Description of a code}

All the code is available on Github:
\begin{verbatim}
https://github.com/nurdin-takenov/lgs_ado
\end{verbatim}

Files with \verb|".nb"| extension are executed in Mathematica. Files with \verb|".py"| extension are executed in Python.

All knots are represented as closures of braids. Braids represented in the following format:
\begin{verbatim}{number of strands, {twist 1, twist 2, ...}}\end{verbatim}

Generator $s_k$ is represented as $k$, $s_k^{-1}$ is represented as $-k$, the braid is read from bottom up.

\begin{figure}[h]
\includegraphics[scale=1.0]{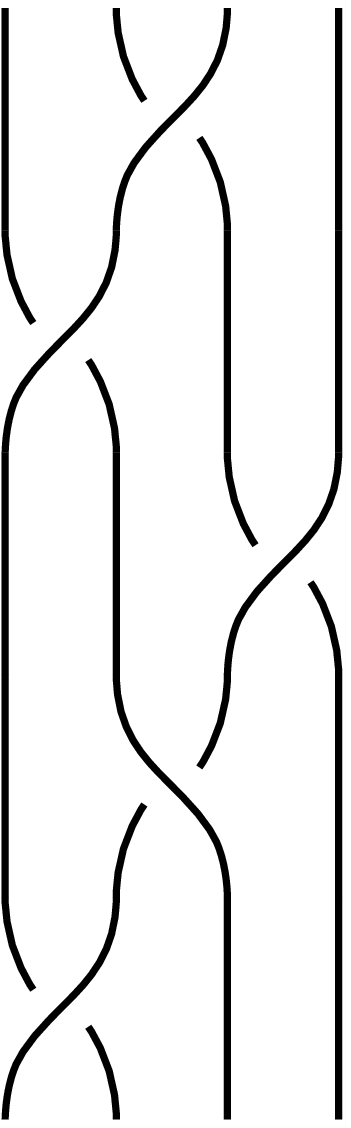}
\caption{The braid \texttt{\{4,\{1,-2,3,1,2\}\}}.}
\end{figure}

Here is the list of the files and their description:

\verb|checking_identities.nb| - this file checks cubic and Ishii relations for $ADO_3$ invariant $R$-matrix \ref{common}.

\verb|s4_elements_generator.py| - this Python code generates elements of $S_4$ \ref{sn}. The results are listed in the file \verb|s4elements.txt|.

\verb|checking_s4.nb| - this file checks the difference between $ADO_3$ and specialized Links-Gould invariant for elements of $S_4$ \ref{sn}.\\

Regarding $S_5$: we need to check 6480 elements, which are broken into 10 types \ref{a5}:

Type 1: $S_4 s_4 s^{-1}_3 s_4$.

Type 2: $S_4 s^{-1}_4 s_3 s^{-1}_2 s_3 s^{-1}_4$.

Type 3: $S_4 s_4 s^{-1}_3 s_2 s^{-1}_3 s_4$.

Type 4: $S_4 s^{-1}_4 w^+ s^{-1}_4$.

Type 5: $S_4 s_4 w^- s_4$.

Type 6: $S_4 s^{-1}_4 w^- s^{-1}_4$.

Type 7: $S_4 s_4 w^+ s_4$.

Type 8: $s_4 w^- s_4 w^- s_4 S_4$.

Type 9: $s_4 w^+ s^{-1}_4 w^+ s_4 S_4$.

Type 10: $s^{-1}_4 w^- s_4 w^- s^{-1}_4 S_4$.

Here $w^+ = s_3 s^{-1}_2 s_1 s^{-1}_2 s_3$, $w^- = s^{-1}_3 s_2 s^{-1}_1 s_2 s^{-1}_3$.

In order to perform calculations faster, the separate code was written for each type of a braid. For each type \verb|k| we have three files:

\verb|calculate_beginning_for_ADO_S5_type_k.nb| 

\verb|calculate_beginning_for_LGS_S5_Type_k.nb|

\verb|checking_S5_type_k.nb|

\begin{figure}[h]
\includegraphics[scale=0.5]{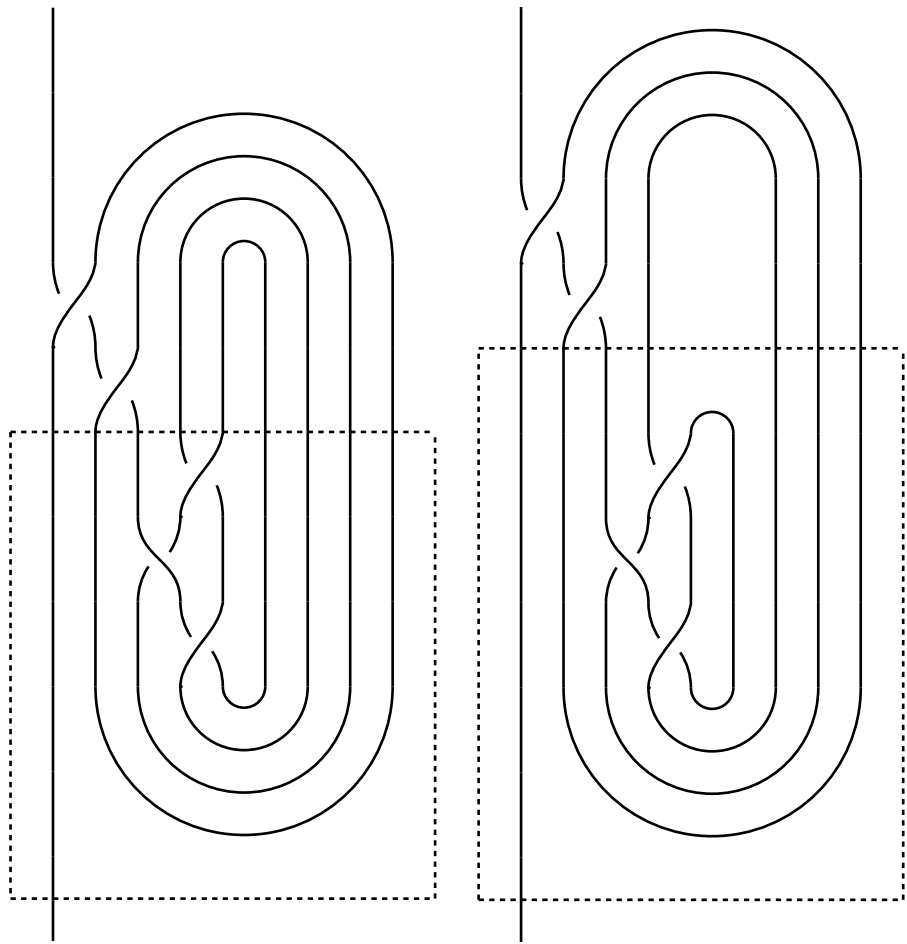}
\caption{The long knot corresponding to the braid \texttt{\{5,\{4,-3,4,2,1\}\}}.}
\end{figure}

Here is the purpose of these files: \verb|checking_S5_type_k.nb| is the main file, that checks equality of specialized Links-Gould and $ADO_3$ invariants for the braids of type \verb|k|. However, it needs two other files. They calculate the operator corresponding to the beginning of the braid.

Let's consider an example. Suppose we want to check braids of the type 1, $S_4 s_4 s^{-1}_3 s_4$ and in particular \verb|{5,{2,1,4,-3,4}}|. Using first Markov move we can change it into braids of the type $s_4 s^{-1}_3 s_4 S_4$, our particular braid will become 
\verb|{5,{4,-3,4,2,1}}|. Then we need to consider a long knot corresponding to this braid. It is shown on the picture on the left. We can isotope the knot and get the picture on the right. The auxiliary files \verb|calculate_beginning_for_ADO_S5_type_k.nb| and \verb|calculate_beginning_for_LGS_S5_Type_k.nb|
 calculate operators corresponding to the tangle encircled by dashed lines. It is the same for the knots of the specific type. This information the is used in the file \verb|checking_S5_type_k.nb| to speed up calculations. Also, since the beginning for knots of the specific type is the same, the common part is not listed. 
 
For example to check whether specialized Links-Gould and $ADO_3$ invariants coincide for the knot corresponding to the braid \verb|{5,{4,-3,4,2,1}}|, we run the command \verb|LGSMADO[{4, {2, 1}}]| in the file \verb|checking_S5_type_k.nb|. We are skipping the part \verb|{4,-3,4}| because it's the common part for the braids of the type 1. The blank file without pre-calculated beginnings is \verb|checking_blank_s5.nb|.
 
 \end{document}